\documentclass{tran-l}
 \usepackage{amsmath}
 \usepackage{amssymb}
 \usepackage[noadjust]{cite}
 \usepackage{color}
 \usepackage{curves}
 \usepackage{enumerate}
 \usepackage{epic}
 \usepackage{epsfig}
 \usepackage{float}
 \usepackage{graphics}
 \usepackage{graphicx}
 \usepackage{latexsym}
 \usepackage{multicol}
 \usepackage{subfigure}
 \newtheorem{thm}{Theorem}[section]
 \newtheorem{hyp}{Hypothesis}

 \newtheorem{conj}[thm]{Conjecture}
 \newtheorem{prop}[thm]{Proposition}
\theoremstyle{remark}
 
 \numberwithin{equation}{section}
\DeclareMathOperator{\A}{A}
 \DeclareMathOperator{\bc}{BC}

 \DeclareMathOperator{\diag}{diag}
 \DeclareMathOperator{\dom}{domain}

 \DeclareMathOperator{\id}{id}

 \DeclareMathOperator{\range}{range}

 \DeclareMathOperator{\real}{Re}

 \DeclareMathOperator{\Span}{Span}
 \newcommand{\BCE}{\bc_e(\GR^2,\GR^m)}
 \newcommand{\BCU}{\bc_u(\GR^2,\GR^m)}

 \newcommand{\GC}{\mathbb{C}}

 \newcommand{\GN}{\mathbb{N}}
 
 \newcommand{\GQ}{\mathbb{Q}}
 \newcommand{\GR}{\mathbb{R}}
 \newcommand{\GS}{\mathbb{S}}
 \newcommand{\GSE}{\mathbb{SE}}
 \newcommand{\GSO}{\mathbb{SO}}
 \newcommand{\GT}{\mathbb{T}}
 \newcommand{\GZ}{\mathbb{Z}}

 \newcommand{\omp}{\omega_{\tiny\mbox{per}}}
 \newcommand{\omr}{\omega_{\tiny\mbox{rot}}}
 \newcommand{\ov}{\overline}

 \newcommand{\sdp}{\dot{+}}

 \newcommand{\newl}{\newline\newline}
\begin{document}

\title[Spiral waves and the dynamical system approach]
 {Spiral waves and the dynamical system approach}

\author{ Patrick Boily }

\address{Department of Mathematics and Statistics, University of Ottawa, Ottawa, Ontario, K1N~6N5,
 Canada}

\email{pboily@uottawa.ca}

\thanks{}

%\subjclass{37G40}


\date{}

\dedicatory{}

\commby{}

%%% ----------------------------------------------------------------------

\begin{abstract}
Spirals are common in Nature: the snail's shell and the ordering of seeds in
the sunflower are amongst the most widely-known occurrences. While these are
static, dynamic spirals can also be observed in excitable systems such as heart
tissue, retina, certain chemical reactions, slime mold aggregates, flame
fronts, etc. The images associated with these spirals are often breathtaking,
but spirals have also been linked to cardiac arrhythmias, a potentially fatal
heart ailment.
\par
In the literature, very specific models depending on the excitable system of
interest are used to explain the observed behaviour of spirals (such as
anchoring or drifting). Barkley \cite{B1} first noticed that the Euclidean
symmetry of these models, and not the model itself, is responsible for the
observed behaviour. But in experiments, the physical domain is never Euclidean.
The heart, for instance, is finite, anisotropic and littered with
inhomogeneities. To capture this loss of symmetry, LeBlanc and Wulff
\cite{LW,LeB} introduced forced Euclidean symmetry-breaking (FESB) in the
analysis.
\par To accurately model the physical situation, two basic types of symmetry-breaking perturbations
are used: translational symmetry-breaking (TSB) and rotational
sym\-metry-breaking (RSB) terms.  In this paper, we provide an overview of
currently know results about spiral wave dynamics under FESB.
\end{abstract}

%%% ----------------------------------------------------------------------
\maketitle
%%% ----------------------------------------------------------------------
\section{Introduction}
The spiral is an integral part of Nature: it can be seen in a snail's shell, in
the layout of a sunflower's seeds and in the path of a falcon on a hunt, to
name but a few. These particular instances are fixed in space, but spirals can
also evolve in time: hurricanes and galaxies are common examples that come to
mind. It is, however, rather arduous to conduct experiments on the latter
physical objects, for obvious reasons.
\par On a smaller scale, where experiments are easier to control, spirals
in evolution have also been observed in \textit{excitable media} such as heart
tissue, slime-mold aggregates, the retina or certain chemical reactions (such
as the famed Belousov-Zhabotinsky (BZ) reaction). In these systems, waves
propagate by `exciting' a `cell', which in turn `excites' some of its
neighbours before falling into a `refractory' or `unexcitable' state, followed
by a `resting' state, ready to be `excited' should the wave come its way again.
\par These systems give rise to beautiful images (as can be attested to in
\begin{figure}[t]
\begin{center}\includegraphics[width=165pt]{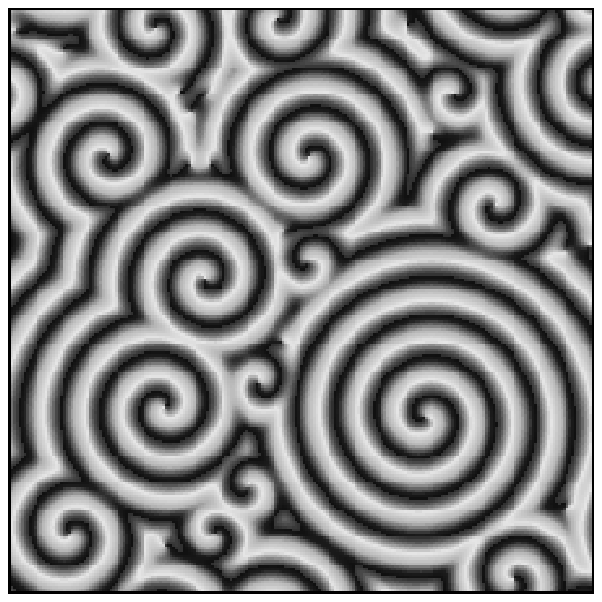}\qquad \includegraphics[width=165pt]{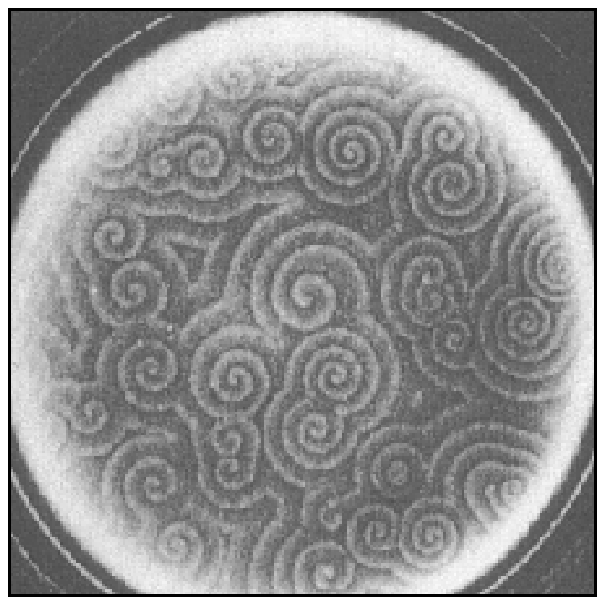}
\caption[Spirals in excitable media.]{\small Spirals in excitable media. On the
left, spirals in a solution of the Ginzburg-Landau equation \cite{Hendrey}, on the right,
spirals in colonies of the slime-mold aggregate \textit{Dictyostelium discoideum} \cite{Ball}.\label{Ballill}}
\end{center}\hrule\end{figure}
 figure~\thefigure). While this in itself might yield enough interest to study them, there is also
 (at least) one serious reason to do so: spiral waves have been linked to
 \textit{cardiac arrhythmias}, \textit{i.e.} to disruptions of
 the heart's normal electrical cycle \cite{W,Wetal}. Most arrhythmias are harmless but if they are
 `re-entrant in nature and [...] occur [in the ventricles] because of the spatial distribution of cardiac
 tissue \cite[p.~401]{KS}', they can seriously hamper the pumping mechanism of the heart and so lead to death.
 As a result, a full understanding of spiral wave dynamics in these media becomes imperative.

%%% ----------------------------------------------------------------------

\section{Historical Perspective}
\label{HP} Numerous experiments and simulations have been performed with
excitable media, see for instance \cite{B1,B2,BKT,Detal,LOPS} (a selected
bibliography can be found in appendix~A). The various `spiral' motions that are
observed are classified according to their \textit{tip path}, an arbitrary
point on the wave front that is followed in time, as can be seen, for instance,
in \addtocounter{figure}{1}figure~\thefigure\addtocounter{figure}{-1}. Some of
the standard possibilities are shown in
\addtocounter{figure}{2}figure~\thefigure.\addtocounter{figure}{-2}
\par In the literature (see \cite{BK,R1,W1} for instance), specific systems of partial differential equations (PDE) have
been used to attempt to explain the observed phenomena: for instance, the
FitzHugh-Nagumo equations \begin{align*}
u_t&=\textstyle{\frac{1}{\varsigma}}(u-\textstyle{\frac{1}{3}}u^3-v)+\Delta u\\
v_t&=\varsigma(u+\beta-\gamma v)
\end{align*} of cardiology, where $\varsigma$,
$\beta$ and $\gamma$ are model parameters, $u$ represents an electric potential
and $v$ a measure of permeability, or the Oregonator
\begin{align*}
u_t&=\textstyle{\frac{1}{\varsigma}}(u-u^2-fv\textstyle{\frac{u-q}{u+q}})+\Delta u\\
v_t&=u-v+D_v\Delta v
\end{align*} for the Belousov-Zhabotinsky reaction, where $f$,
$\varsigma$ and $q$ (small) are the model parameters, $D_v$ is a diffusion
coefficient and $u$ and $v$ represent the concentrations of certain chemical
reactants.\par These two systems of partial differential equations are instances of a
general class of PDE, given by
\begin{figure}[t]
\begin{center}\includegraphics[width=150pt]{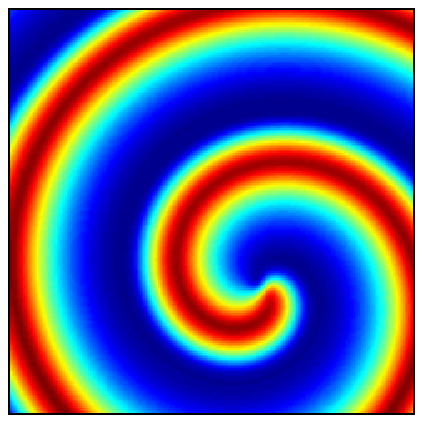}\qquad \includegraphics[width=150pt]{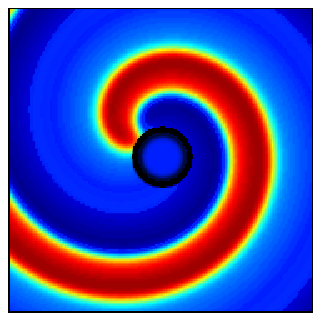}
\caption[Isolated spiral and spiral tip path in excitable
media.]{\small On the left, an isolated spiral in a
RDS; the excited wave front is shown in red.
On the right, the corresponding tip path.}
\end{center}\end{figure}
 \begin{figure}[t]
\begin{center}\includegraphics[width=320pt]{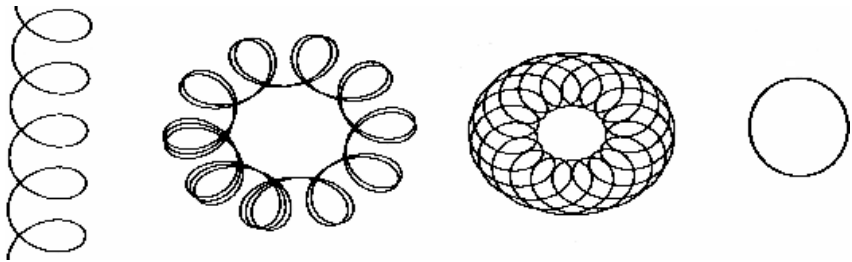}%\includegraphics[width=180pt]{fig2}
\caption[Possible paths of the spiral wave tip.]{\small Possible paths of the spiral wave tip.
From left to right: linear drifting, outer epicycle motion, inner epicycle motion and rigid rotation \cite{LW}.\label{LW111}}\end{center}\hrule
\end{figure}
\begin{equation}u_t(x,t)=f(u(x,t))+D\Delta u(x,t),\label{thefirst}\end{equation} where $x\in \GR^2$, $u:\GR^2\times \GR^+_0\to \GR^m$
is bounded and uniformly continuous, $D$ is an $m\times m$ diagonal matrix and
$f:\GR^m\to \GR^m$ is some sufficiently smooth function. General systems of the form (\ref{thefirst}) are called
\textit{reaction-diffusion system} (RDS).\footnote{Winfree provides a very
complete survey of their use as models \cite{W1}.}\newpage\noindent Four types of
solutions of (\ref{thefirst}) are of particular interest in the context of spiral wave dynamics:
\begin{description}
\item[Rotating waves](RW) are rigidly rotating periodic solutions that are fixed in a
\textit{co-rotating} frame of reference (\textit{i.e.}
the frame rotates uniformly with the same frequency as the solution). In
physical and numerical experiments, the tip path is circular. RW are sometimes
called \textit{vortices} or \textit{rotors} in the literature
\cite{W1,Wetal}.
\item[Traveling waves](TW) are
linearly propagating solutions that are fixed in a \textit{co-translating} frame of
reference (\textit{i.e.} the frame translates
linearly and uniformly with the solution). In experiments, the tip path of such
a solution is a line. Strictly speaking, TW (or \textit{retracting tip waves} \cite{BK}) are not spiral
waves as they do not have a rotating component.
\item[Modulated rotating waves] (MRW) are two-frequency quasi-periodic
solutions that are periodic in a co-rotating frame of
reference that rotates uniformly with one of the
frequencies of the solution. The tip path of such a solution is a closed
epicycle when the ratio of the frequencies is rational; otherwise the tip path
densely fills a ring over time, with an epicycle-like motion. In the
literature, MRW are sometimes called \textit{meandering waves}.
\item[Modulated traveling waves] (MTW)
are rotating solutions, superimposed with a linearly propagating motion, that
are periodic in a co-translating frame of reference
that travels uniformly with the linear component of the solution. The tip path
of such a solution is a helix-shaped two-dimensional curve.
\end{description}
As an example of a RDS in which these occur, consider Barkley's
system:
\begin{align}\label{pd}
\begin{split}
u_t&=\textstyle{\frac{1}{\varsigma}u(1-u)\left(u-\frac{v+b}{a}\right)+\Delta u}\\
v_t&=u-v,
\end{split}
\end{align} where $a,b$ and $\varsigma$ are system parameters with $\varepsilon$ small \cite{B2}.
\begin{figure}[t]
\begin{center}
\includegraphics[width=225pt]{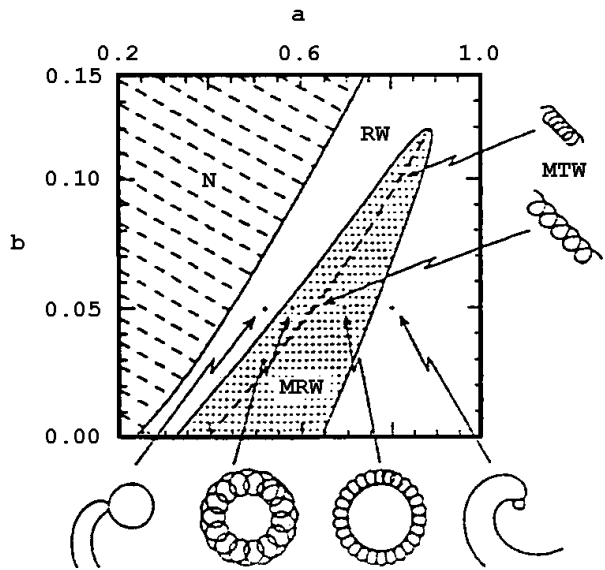}
\caption[Bifurcation diagram of spiral wave dynamics in (\ref{pd}).]{\small
\label{thingy} Bifurcation diagram of spiral wave dynamics in (\ref{pd}) \cite{B2}.}\label{Bark1}
\end{center}\hrule\end{figure} Figure~\thefigure \ shows a bifurcation diagram
of the spiral dynamics of (\ref{pd}) for $\varepsilon=\frac{1}{50}$. There are
three regions of interest labeled \texttt{N}, \texttt{RW} and \texttt{MRW}.
Note that this last region is divided in two sub-regions by a curve labeled
\texttt{MTW}. \par In \texttt{N}, no wave propagation is observed; in
\texttt{RW}, observed solutions are RW and in \texttt{MRW}, observed solutions
are MRW, with petality\footnote{The orientation of the
spiral ``petals''.} determined by the side of the curve \texttt{MTW} on which
the parameters fall. The intersection of this curve with the boundary of
\texttt{MRW} is a point that deserves special consideration: in every one of
its neighbourhoods, the three basic types of spiral behaviours can be seen.
\par An \textit{a priori} surprising feature of reaction-diffusion systems is
that figure~\thefigure \ is a generic bifurcation diagram: most experimental results are strikingly similar \cite{BK, W1}; this
suggests they are in fact a consequence of excitable media and
their geometry, and not of the particular models that aim to describe the
dynamics \cite{B1,GLM,LeB,LW}.\footnote{ It should be noted that
reaction-diffusion systems are not the sole models of excitable media, nor were
they the first: Wiener and Rosenblueth originally defined and modeled excitable
media using cellular automata \cite{WR}.} This is where the dynamical system approach enters the picture.

%%% ----------------------------------------------------------------------

\section{Equivariant Vector Fields and Reaction-Diffusion Systems}
\label{Wulffsection} Let $\Gamma$ be a group acting linearly on a vector space
$X$. A function $f:X\to X$ is \textit{$\Gamma-$equivariant} if it commutes with
the action of $\Gamma$, \textit{i.e.} $$\gamma\cdot f(x)=f(\gamma\cdot x),\quad
\forall \gamma\in \Gamma,x\in X.$$ Equivariant vector fields (with compact
$\Gamma$) have been studied by many authors: notable amongst them are
Golubitsky and Schaeffer \cite{GS}, Golubitsky, Stewart and Schaeffer
\cite{GSS} and Vanderbauwhede \cite{V,V1}. The main feature of these
$\Gamma-$equivariant vector fields is that whenever $x(t)$ is a solution of
$\dot{x}=f(x)$, so is $\gamma x(t)$, for all $\gamma\in \Gamma$.\par The
\textit{special Euclidean group} $\GSE(2)=\GC\sdp \GS^1$ is a non-compact
subset of all the distance-preserving transformations of the plane, with
multiplication defined by
\begin{equation}\label{SEmult} (p_1,\varphi_1)\cdot (p_2,\varphi_2)=(e^{i\varphi_1}p_2+p_1,\varphi_1+\varphi_2),
\quad \forall (p_1,\varphi_1),(p_2,\varphi_2)\in \GSE(2).\end{equation} It acts
on the space of bounded uniformly continuous functions from $\GR^2$ to $\GR^m$,
which we will denote by $\BCU$, according to
\begin{equation}\label{thesecond}(\gamma \cdot
v)(x)=((p,\varphi)\cdot v)(x)=v(R_{-\varphi}(x-p)), \quad \forall
(p,\varphi)\in \GSE(2),\end{equation} where $R_{\theta}$ represents a rotation
by angle $\theta$ around the origin. Reaction-diffusion systems on $\BCU$ are
$\GSE(2)-$equivariant under the action of (\ref{thesecond}), but that action is
not smooth over $\BCU$: the problem arises with rotations, as a small shift in
$\theta$ produces a large displacement at far distances \cite{Wulff}. However,
there is a closed set $\BCE\subsetneq \BCU$ over which (\ref{thesecond}) is
smooth \cite{Wulff}.

\subsection{Abstract Differential Equations} In order to determine  $\BCE$,
Wulff uses the following RDS paradigm \cite{W}. Consider
\begin{equation}\label{rdcmrt} u_t(x,t)=\tilde{D}\Delta
u(x,t)+f(u(x,t),\varsigma),\end{equation} where $x\in \GR^2$, $u:\GR^2\times
\GR^+_0\to \GR^m$, $\tilde{D}\geq 0$ is a diagonal matrix, $\varsigma\in
\GR^M$, $\Delta$ is the Laplacian and $f$ is $C^{k+2}$ for some $0 \leq k\leq
\infty$. If $\det\tilde{D}\neq 0$, let $Y=\BCU$ be the Banach space of
\textit{uniformly continuous, bounded functions from $\GR^2$ to $\GR^m$}.
Otherwise, as long as $f$ satisfies some additional growth conditions, the
choice $Y=L^2(\GR^2,\GR^m)$ can also be used, with slight variations (see
\cite{SSW2,H} for details).\footnote{However, physical considerations demand
that spiral waves in an infinitely extended medium be located in $\BCU$
\cite{Wulff}.} \par The \textit{semi-linear differential equation on $Y$
associated to $(\ref{rdcmrt})$} is the abstract differential equation
\begin{equation}\label{slde} \frac{du}{dt}=-Au+F(u,t,\varsigma),\end{equation} where $F(u,t,\varsigma)=f(u(\cdot,t),\varsigma)$
and \begin{equation}\label{theAtheF}A=\diag(-d_1\Delta,\ldots,
-d_m\Delta).\end{equation} Solutions of (\ref{slde}) are in one-to-one
correspondence with solutions of (\ref{rdcmrt}).\par In the remainder of this
section, we assume the reader is familiar with basic definitions and results
from operator theory and abstract differential equations (see
\cite{H,Wulff,AKVA} for details and definitions).
\begin{prop}\label{sectprop} \textsc{\cite{Wulff}} Let $A$ be given by $(\ref{theAtheF})$, $Y=\BCU$,
$\alpha\in (\frac{1}{2},1)$ and $A_1=\id_Y-\tilde{D}\Delta$, where $\tilde{D}$
is as in $(\ref{rdcmrt})$. Then $A$ is sectorial in $Y^{\alpha}$  and
$\frac{\partial}{\partial x_1}A_1^{-\alpha},\frac{\partial}{\partial
x_2}A_1^{-\alpha}$ are bounded on $Y$. \end{prop}
\subsection{Existence and Uniqueness of
Solutions} That (\ref{slde}) has classical solutions is shown by the following
theorem.
\begin{thm}\label{exun} \textsc{\cite{H,Wulff}}
Let $Y$, $A$, $A_1$ and $\alpha$ be as in proposition $(\ref{sectprop})$, $U$
be a subset of $Y^{\alpha}\times\GR\times \GR^M$ and $F$ be as in
$(\ref{slde})$, locally Lipschitz in its first variable and continuous in the
remaining variables. Then, for any $(u_0,t_0,\varsigma)\in U$, $(\ref{slde})$
has a unique classical $C^{k+2}-$sol\-ution $u(t;u_0,t_0,\varsigma)$ on
$[t_0,t_1]$, where $t_1=t_1(u_0,t_0,\varsigma)>t_0$ and $0\leq k\leq\infty$ is
the smoothness of the nonlinearity $f$ in $(\ref{rdcmrt})$.
\end{thm}
\subsection{$\GSE(2)-$Equivariance of Solutions} Wulff then shows that these solutions
have a $\GSE(2)-$equivariant structure. A \textit{smooth local semi-flow}
$\{\Phi_{t}\}_{t\geq 0}$ on a Banach space $X$ is a smooth family of operators
satisfying $\Phi_{0}=\id_{X}$ and
\begin{equation}\label{semiflow}\Phi_{t+s}=\Phi_{t,\varsigma}\circ
\Phi_{s}=\Phi_{s}\circ \Phi_{t}\quad\mbox{for all }s,t\geq 0.\end{equation} If
furthermore $\Phi_{t}x\to x$ as $t\to 0^+$ for all $x\in X$, then
$\{\Phi_t\}_{t\geq 0}$ is a $C^0-$\textit{semi-group} on $X$. A \textit{smooth
semi-group} on $X$ is a $C^0-$semi-group for which the map
$\Psi_{x}:(0,\infty)\to X$ defined by $\Psi_{x}(t)=\Phi_{t}(x)$ is smooth for
all $x\in X$.
\par The \textit{infinitesimal generator $L_{T}$} of a smooth semi-group
$\{T_t\}_{t\geq 0}$ on a Banach space $X$ is
\begin{equation} \label{theinfinitesimal} L_{T}x=\lim_{t\to 0^+}\frac{1}{t}\left(\Phi_{t}x-x\right),\end{equation}
whenever the limit exists. \par The \textit{special Euclidean group} $\GSE(2)$
plays an important role in the theory of spiral waves. For now, we assume it is
parameterized as $\GSE(2)=\GSO(2)\sdp \GR^2,$ with multiplication given by
\begin{equation}\label{classs}(R_1,S_1)\cdot
(R_2,S_2)=(R_1R_2,S_1+R_2S_2).\end{equation} The \textit{standard
$\GSE(2)-$action on $\BCU$} is \begin{equation}\label{GSEgacmrt}
(R,S)\left(v(x,t)\right)=v\left(R^{-1}(x-S),t\right).\end{equation} Note that
$\GSE(2)$ is generated by the families $\{S^1_{\mu}\}$, $\{S^2_{\nu}\}$
and $\{R_{\theta}\}$, where $$S^1_{\mu}=\begin{pmatrix} \mu \\
0\end{pmatrix}, \quad S^2_{\nu}=\begin{pmatrix} 0 \\ \nu\end{pmatrix}
\quad\mbox{and} \quad R_{\theta}=\begin{pmatrix} \cos \theta & -\sin \theta \\
\sin \theta & \cos\theta \end{pmatrix}$$ satisfy (\ref{semiflow}). Of these,
only $\{S^1_{\mu}\}_{\mu\geq 0}$ and $\{S^2_{\nu}\}_{\nu\geq 0}$ are smooth
semi-groups on $\BCU$; their respective infinitesimal generators are
\begin{equation*}L_{S^1}=-\frac{\partial}{\partial
x_1}\quad\mbox{and}\quad L_{S^2}=-\frac{\partial}{\partial x_2}.
\end{equation*} On the other hand,
$\{R_{\theta}\}$ is not a smooth semi-group (see \cite[lemma 2.14]{Wulff}); the
action of $\GSE(2)$ on $\BCU$ is not even continuous. Thus, $\BCU$ is not a
suitable space over which to define (\ref{rdcmrt}). \par This obstacle is
overcome as follows. The formal evaluation of (\ref{theinfinitesimal}) yields
\begin{equation*}L_R=x_2\frac{\partial}{\partial
x_1}-x_1\frac{\partial}{\partial x_2}. \end{equation*} Now let $\tilde{Y}=\BCE$
be the topological closure of $\dom(L_R)$ in $\BCU$. The reaction-diffusion
system (\ref{rdcmrt}) is well-posed on $\tilde{Y}$, and proposition
(\ref{sectprop}) and theorem \ref{exun} still hold after substituting
$\tilde{Y}$ for $Y$ \cite{Wulff}. Furthermore, (\ref{GSEgacmrt}) is continuous
on $\tilde{Y}$ and the following result holds.
\begin{thm} \textsc{\cite{DMcK,T,Wulff}} \label{theac} The semi-flow $\Phi_{t,\varsigma}$ generated by $(\ref{rdcmrt})$
commutes with the restricted action $(\ref{GSEgacmrt})$ of $\GSE(2)$ over
$\BCE$.
\end{thm}
In particular, the other results of her thesis hold as long as
$\Phi_{t,\varsigma}$ is a smooth $\GSE(2)-$equivariant semi-flow (not
necessarily generated by a RDS) on some suitable  Banach space.
\section{Barkley's Insight}
\label{BI} Barkley was the first to realize that the Euclidean symmetry
discussed in the previous section (as opposed to the specifics of a given
model) could explain most of the spiral wave dynamics observed in experiments
and simulations, and succinctly presented in figure~\thefigure\  \cite{B1,B2}.
His key observation rests on the fact that for any reaction-diffusion system,
the linearization at a RW at the onset of a Hopf bifurcation (hence at the
boundary of \texttt{MRW}) has five isolated leading eigenvalues on the
imaginary axis: $\lambda_R=0$ (due to rotational symmetry), $\lambda_T=\pm
i\omega$ (due to translational symmetry) and $\lambda_B=\pm i\beta_0$
(responsible for the Hopf bifurcation from RW to MRW or \textit{vice-versa}).
\par The pairs of complex conjugate eigenvalues $\lambda_T$ and $\lambda_B$
can be made to coincide by varying two or more system parameters. The
corresponding (interesting) codimension-two point is then found to lie
precisely at the intersection of \texttt{MTW} and the boundary of \texttt{MRW}.
\subsection{Linear Stability Analysis at a RW}
In his ground-breaking paper \cite{B1}, Barkley considers the Oregonator-like
system
\begin{align*}
\frac{\partial u_1}{\partial t}&=\Delta u_1+\frac{1}{\varepsilon}u_1(1-u_1)\left(u_1-\frac{u_2+b}{a}\right)\\
\frac{\partial u_2}{\partial t}&=\delta \Delta u_2+u_1-u_2,
\end{align*} where $a,b,\varepsilon$ are parameters with $\varepsilon$ small and $\delta\in [0,1]$ constant.
In vector form, this system may be written as
\begin{align}\label{Barksys}
\frac{\partial u}{\partial t}=\hat{\delta} \Delta u +f(u),
\end{align} where $u= (u_1,u_2)^{\!\top\!},$ $\hat{\delta}=\diag (1,\delta)$ and $f$ contains the remaining terms.
The boundary conditions $\partial _r u=0$ is taken on a circle of radius $R>0$.
With this set-up, (\ref{Barksys}) is $\GSE(2)-$equivariant under the action of
(\ref{GSEgacmrt}), on some suitable Banach function space, as $R\to \infty$.
\par To find RW solutions of (\ref{Barksys}), \textsl{i.e.} solutions for which
$(\partial_t+\omega
\partial_{\theta})u\equiv 0$\index{rotating wave} for some speed of rotation
$\omega$, it suffices to solve the eigenvalue problem
\begin{align}\begin{split}\label{Barkeig} F(u)&= 0 \\
DF(u)\tilde{u}&=\lambda \tilde{u},
\end{split}\end{align} where $F(u)=\hat{\delta}\Delta u +\omega \partial_{\theta}u+f(u)$ and
$DF(u)= \hat{\delta}\Delta+\omega\partial_{\theta}+Df(u)$. Any $\lambda$
solving (\ref{Barkeig}) corresponds to an eigenvalue of the linearization of
(\ref{Barksys}) at the RW solution. \par Using fast and efficient numerical
methods, Barkley shows that three of the five leading
eigenvalues\footnote{Eigenvalues with largest real part.} lie on the imaginary
axis. Indeed, the rotational symmetry of (\ref{Barksys}) forces $\lambda_R=0$
(with corresponding eigenmode $\tilde{u}_R=\partial_{\theta}u$, where $u$ is
the spiral solution of the first equation in (\ref{Barkeig})); the
translational symmetry of (\ref{Barksys}) imposes $\lambda_T=\pm i\omega$ (with
corresponding eigenmode $\tilde{u}_T=\partial_{x}u\pm i\partial_y u$, where $u$
is as above). Note that this holds in spite of the fact that the boundary
condition breaks the Euclidean symmetry: for sufficiently large domain, the
real part of $\lambda_T$ is numerically indistinguishable from
zero.\footnote{Barkley provides some very strong estimates to that effect.}\par
There is a last pair of complex conjugate leading eigenvalues
$\lambda_B=\alpha(a)\pm i\beta(a)$ that crosses the imaginary axis for some
prescribed $a=a^*$, leading to a Hopf bifurcation or `spiral wave instability',
in which MRW are observed.
\par All the remaining eigenvalues have negative real part and so do not affect
spiral dynamics. As a result, the five leading eigenvalues are isolated in the
spectrum and so any pair $(u,\omr)$ that solves (\ref{Barkeig}) is \emph{not}
part of a continuum of solutions with continuously varying shapes or speed of
rotation. \par These results are in fact model-independent; as long as the
reaction-diffusion equations governing the field $u$ are $\GSE(2)-$equivariant,
the five leading eigenvalues will have the above properties.
\subsection{The Ad-Hoc Model}

Based on these observation, Barkley \cite{BK} constructed an \textit{ad hoc}
$5-$dimensional system of ordinary differential equations (ODE) which
replicates the above resonant Hopf bifurcation:
\begin{align}\label{BarKsys}
\begin{split}
\dot{p}&=v \\
\dot{v}&= v\left[f(|v|^2,w^2)+iwh(|v|^2,w^2)\right] \\
\dot{w}&=wg(|v|^2,w^2)
\end{split}
\end{align}
where $p,v\in \GC$, $w\in \GR$ and
$$f(\xi,\zeta)=-\frac{1}{4}+\alpha_1\xi+\alpha_2\zeta-\xi^2,\quad
g(\xi,\zeta)=\xi-\zeta-1,\quad \mbox{and}\quad h(\xi,\zeta)=\gamma_0$$ for some
$\gamma_0\in \GR$. The variable $p$ represents the position of the spiral tip, while $v$ is its
linear velocity and $\gamma_0w$ its instantaneous rotational frequency rate. This system has RW solutions that undergo a Hopf bifurcation
to MRW solutions and it
also has a codimension-two resonant Hopf point. Furthermore, it is
equivariant under the distance-preserving planar transformations generated by
$$\textstyle{R_{\gamma}\begin{pmatrix}p \\ v \\ w\end{pmatrix}=\begin{pmatrix}e^{i\gamma} p \\ e^{i\gamma}v \\ w\end{pmatrix}
\quad\mbox{and}\quad T_{\alpha,\beta}\begin{pmatrix}p \\ v \\
w\end{pmatrix}=\begin{pmatrix}p+\alpha+i\beta \\ v \\ w\end{pmatrix}}, $$ where
$R_{\gamma}$ and $T_{\alpha,\beta}$ represent respectively a rotation by angle
$\gamma$ around the origin and a
translation by the vector $\alpha+i\beta$.
\par Note the absence of $p$ in the right-hand side of (\ref{BarKsys}) as position
plays no role in Euclidean systems.
\begin{figure}[t]
\begin{center}
\includegraphics[width=210pt]{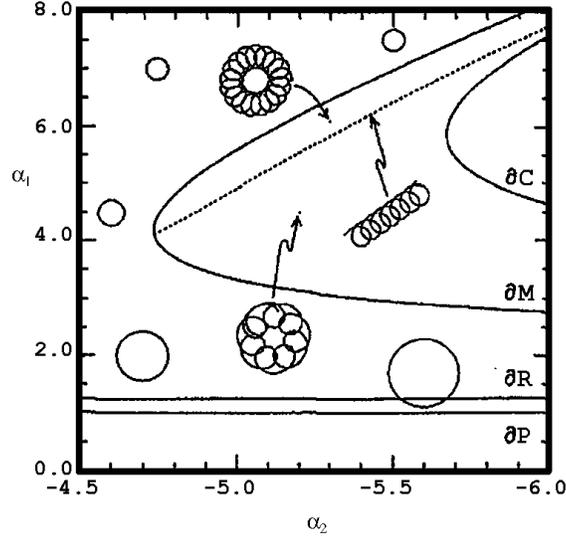}
\caption[Bifurcation diagram of (\ref{BarKsys}).]{\small Bifurcation diagram of
(\ref{BarKsys}) \cite{BK}.}\end{center}\label{BKfig}\hrule
\end{figure} Figure~\thefigure\ shows the bifurcation diagram of (\ref{BarKsys}) for $\gamma_0=5.6$.
The similarities with
\setcounter{figure}{4}figure~\thefigure\setcounter{figure}{5} are readily
apparent, in particular when it comes to the presence of a codimension-two
organizing center around which the three types of spirals can be found.
\section{The Dynamical System Approach}
 Then, in what has been hailed a ``major mathematical work on spirals
\cite{GLM}'', Wulff~\cite{Wulff} rigourously proved that the resonant unbounded
growth observed by many authors (such as \cite{B1,BK}) does indeed occur near
the codimension-two point. \newpage\noindent The following result, the center
manifold reduction theorem of Sandstede, Scheel and Wulff
\cite{SSW1,SSW2,SSW3,FSSW}, remains, in the author's opinion, both the most
technical general results on spiral wave dynamics  and its most fruitful ally
in applications. It is an extension to non-compact symmetry group of Krupa's
\cite{KR} center bundle construction for relative equilibria and periodic
solutions.
\subsection{The Center Manifold Reduction Theorem (CMRT)}
The CMRT helps provide a rigorous link between spiral solutions of
(\ref{rdcmrt}) and Barkley's \textit{ad hoc} model~(\ref{BarKsys}). Set
$\mathfrak{se}(2)=\mathfrak{so}(2)\times \GR^2$, where $\mathfrak{so}(2)$ is
the Lie algebra of $\GSO(2)$, consisting of the $2\times 2$ anti-symmetric
matrices on $\GR$. As a one-dimensional vector space,
$$\mathfrak{so}(2)=\Span\left\{\begin{pmatrix}0 & 1
\\ -1 & 0\end{pmatrix}\right\}=\Span \{J_2\}.$$ Define $\exp_{\mathfrak{so}(2)}:\mathfrak{so}(2)\to
\GSO(2)$ by
$$\exp_{\mathfrak{so}(2)}\left(bJ_2\right)=\begin{pmatrix}\cos b & \sin b \\
-\sin b & \cos b\end{pmatrix}.$$ Then, $\mathfrak{se}(2)$ is the Lie algebra of
$\GSE(2)$, when endowed with commutator and exponential maps as defined by
(\ref{commumap}) and (\ref{expumap}) below. Let $I_2$ be the $2\times 2$
identity matrix. Then
\begin{align}\label{commumap}
[(r_1,s_1),(r_2,s_2)]&=(r_1r_2-r_2r_1,r_1s_2-r_2s_1) \\
\exp((r,s)t)&=(\exp_{\mathfrak{so}(2)}(rt),r^{-1}(\exp_{\mathfrak{so}(2)}(rt)-I_2)s),\label{expumap}
\end{align} for $t\in \GR$, $r_j\in \mathfrak{so}(2)$ and $s_j\in \GR^2$, $j=\varnothing, 1, 2$ \cite{FSSW,SSW2}.
\par Now, consider a RDS of the form (\ref{rdcmrt}). For a fixed $\varsigma$, a
\textit{relative equilibrium} of (\ref{rdcmrt}) is a solution $u(x,t)$ whose
time orbit, or semi-flow orbit, is contained in its group orbit under
$\GSE(2)$. More precisely, it satisfies
\begin{equation}\label{releqcmrt}u(x,t)=\exp((r,s)t)u(x,0)\end{equation}
for some $(r,s)\in \mathfrak{se}(2)$.  The \textit{isotropy subgroup} of such
solutions is
$$\Sigma_{u}=\{\sigma\in \GSE(2):\sigma u(x,t)=u(x,t)\}.$$ If
$\Sigma_{u}\cong \GZ_{\ell}$, then $u(x,t)$ is an $\ell-$\textit{armed spiral}.
\par According to the definitions of section \ref{HP}, RW are relative
equilibria of (\ref{rdcmrt}) with no translation component; after a change of
coordinates bringing the center of rotation to the origin, (\ref{releqcmrt})
becomes
\begin{equation}\label{releqcmrtrrw}u_{*}(x,t)=(\exp_{\mathfrak{so}(2)}(r_{*}t),0)u_{*}(x,0)\end{equation}
for some non-trivial $r_{*}\in \mathfrak{so}(2)$. \par Similarly, a
\textit{relative periodic solution} of (\ref{rdcmrt}) is a solution
$$u(x,t)=\exp((r,s)t)w(x,t),$$ for some $(r,s)\in \mathfrak{se}(2)$, where $w$ is
a non-constant $T-$periodic function in $t$. According to the definitions of
section \ref{HP}, MRW are relative periodic solutions of (\ref{rdcmrt}) with no
translation component; after a change of coordinates bringing the center of
rotation to the origin, (\ref{releqcmrt}) becomes
\begin{equation}\label{relpersolcmrtmrw} u^{*}(x,t)=(\exp_{\mathfrak{so}(2)}(r^{*}t),0)w(x,t)\end{equation} for some non-trivial
$r^{*}\in \mathfrak{so}(2)$ and a non-constant $T-$periodic function
$w$.\footnote{While we focus mainly on RW and MRW, the theorems of this section
can easily be adapted for TW and MTW.} The isotropy subgroup of $u^*$ is
defined as for RW, and the interpretation is identical. \par Standard center
bundle results (such as those presented by Krupa \cite{KR}) cannot be applied
to (\ref{rdcmrt}) because $\GSE(2)$ is not compact (a small rotation at the
origin will produce a large displacement away from the origin). However, it can
be shown that relative equilibria and relative periodic solutions are members
of $\BCE$, and so, as was seen in section \ref{Wulffsection}, that the action
of $\GSE(2)$ defined by (\ref{GSEgacmrt}) is continuous on RW and MRW. \newl
The two hypotheses that follow allow the resolution of some technical
difficulties within the CMRT.
\begin{hyp}\label{existence} For the parameter $\varsigma_*$ $[$resp.
$\varsigma^*$$]$, assume $u_*$ $[$resp. $u^*$$]$ is a $1-$armed (normally
hyperbolic) RW $[$resp. MRW$]$ of $(\ref{rdcmrt})$ with $0\neq r_*$ as in
$(\ref{releqcmrtrrw})$ $[$resp. with $0\neq r^*$ and $w$ as in
$(\ref{relpersolcmrtmrw})$$]$. If $\tilde{D}$ is singular, assume further that
$u_*$ $[$resp. $u^*$$]$ is $k+2-$times uniformly continuously differentiable.
\end{hyp}
Scheel \cite{S} has shown that such rotating wave solutions can arise from Hopf
bifurcations in a large class of planar reaction-diffusion equations.
\begin{hyp}\label{spectral} Assume that $\{\mu:|\mu|\geq 1\}$ is a spectral set for
the linearization $\exp_{\mathfrak{so}(2)}(-r_*)D\Phi_{1,\varsigma_*}(u_*)$
$[$resp. $\exp_{\mathfrak{so}(2)}(-r^*)D\Phi_{1,\varsigma^*}(u^*)$$]$ and that
$$\dim(\range(P_*))=3,$$ $[$resp. $\dim(\range(P^*))=5$$]$, where $P_*$ $[$resp.
$P^*$$]$ is the spectral projection associated to $u_*$ $[$resp. $u^*$ $]$.
\end{hyp}
That this second hypothesis \textit{can} hold has been verified numerically by
Barkley~\cite{B1} (see section \ref{BI}). It should be noted, however, that
Scheel \cite{S} has also shown that this hypothesis fails to hold for a large
class of asymptotically Archimedean spiral waves. \newpage\noindent The
following result establishes the existence of an invariant center manifold that
is contained in an $\GSE(2)-$invariant neighbourhood of the group orbit of
$u_*$ [resp.~$u^*$].
\begin{thm}\label{thmmm1} \textsc{(}\textsc{\cite{SSW2}}, $\mathrm{theorem\ 4,\
p.~142}$\textsc{)} For any $\varsigma$ close enough to $\varsigma_*$ $[$resp.
$\varsigma^*$$]$, there exists an $\GSE(2)-$invariant, locally
semi-flow-invariant manifold $M_{\varsigma}^{\mbox{\tiny cu}}$. Both
$M_{\varsigma}^{\mbox{\tiny cu}}$ and the action of $\GSE(2)$ on
$M_{\varsigma}^{\mbox{\tiny cu}}$ are $C^{k+1}$ and depend $C^{k+1}-$smoothly
on $\varsigma$. Furthermore, $M_{\varsigma}^{\mbox{\tiny cu}}$ contains all
solutions which stay close to the group orbit of $u_*$ $[$resp. $u^*$$]$ for
all negative times. Finally, $M_{\varsigma}^{\mbox{\tiny cu}}$ is locally
exponentially attracting. Furthermore, the manifold $M_{\varsigma}^{\mbox{\tiny
cu}}$ is diffeomorphic to the bundle $V_*=\GC\times \GS^1$ $[$resp.
$V^*=\GC\times \GT$$]$.
\end{thm}

\subsection{The CMRT and its Applications}
In~\cite{GLM}, Golubitsky, LeBlanc and Melbourne describe the structure of the
equations on $M_{\varsigma}^{\mbox{\tiny cu}}$ assuming the spiral waves have
trivial isotropy subgroup. Generally, the essential dynamics for Hopf
bifurcation from $\ell-$armed spirals are analyzed \textit{via} a
$5-$dimensional system of ODE on the center bundle $\GSE(2)\times \GC$
describes. For $\ell=1$, the general system reduces to the center bundle
equations
\begin{align}\label{thethird}
\begin{split}
\dot{p}&=e^{i\varphi}F^p(q,\overline{q})\\
\dot{\varphi}&=F^{\varphi}(q,\overline{q})\\
\dot{q}&=F^q(q,\overline{q}),
\end{split}
\end{align} where $p,q\in \GC$, $\varphi\in \GS^1$, $F^{\varphi}(0)=\omr\in \GR$, $F^q(0)=0$ and
$DF^{q}(0)=i\omp$ is purely imaginary.\footnote{The frequencies $\omr$ and
$\omp$ in (\ref{thethird}) play similar roles to the parameters $\omega$ and
$\beta_0$ in \cite{B1}.} The Euclidean action on $\GSE(2)\times \GC$ is given
by
\begin{equation}\label{thefourth} (x,\theta)\cdot
(p,\varphi,q)=(e^{i\theta}p+x,\varphi+\theta,q),\quad \forall(x,\theta)\in
\GSE(2).\end{equation} The analysis of (\ref{thethird}), the titular
\textit{dynamical system approach}, allows the authors to recover the results
of Barkley and Wulff concerning the Hopf bifurcation from a RW and resonant
growth by considering a parameterized version of (\ref{thethird}); a quick note
on Bogdanov-Takens bifurcation from $1-$armed spirals is also provided. If
$\ell>1$, the structure of the 5 dimensional center bundle equations changes,
but, as an ODE system, it retains $\GSE(2)-$equivariance under the action
\begin{equation}\label{THEaction} (x,\theta)_m\cdot (p,\varphi,q)= (e^{i\theta}p+x,\varphi+\ell \theta,e^{im\theta} q),\quad \forall(x,\theta)\in
\GSE(2),\end{equation} for some fixed $m\in
\{0,\ldots,[\ell/2]\}$.\footnote{This action is consistent with
(\ref{thefourth}).}
\begin{thm} \textsc{(}\textsc{\cite{GLM}}, $\mathrm{section\ 5,\
p.~571}$\textsc{)} Let $\dot{y}=N(y,\mu)$ be the center bundle equations for an
$\ell-$armed spiral, parameterized by $\mu\in \GR$. There is a unique parameter
value $\mu_0$ at which the spiral undergoes a codimension-two bifurcation to
resonant growth if and only if $\ell$ and $m$ are coprime.\end{thm}
\par In a subsequent paper \cite{GLM2}, the authors again
use the dynamical system approach to show that while $\GSO(2)-$symmetry alone
may explain rotating waves, Euclidean symmetry is necessary in order to observe
the unbounded growth of Barkley and Wulff, as well as to explain the full
bifurcation diagram of section \ref{HP}.\par Still, some observed spiral
behaviours (see next section for a partial list) are left unexplained by a
careful analysis of the center bundle equations (\ref{thethird}); this obstacle
is overcome through the introduction of forced Euclidean symmetry-breaking.

\section{The Effects of Forced Euclidean Symmetry-Breaking}

Physical experiments can never be perfectly Euclidean, if only because of their
finite nature. In the heart, this reality is obvious. Cardiac tissue is
\textit{anisotropic} (\textit{i.e.} heart fibres have a preferred orientation
and electrical conductivity is direction-dependent). Furthermore, tissue
distribution is not uniform: there are zones of relatively high density that
affect cardiac activity \cite{KS}.
\par Similarly, introducing light pulses in a light-sensitive BZ
reaction changes the geometry of the system. Moreover, the boundary cannot be
ignored when the size of the spiral core is `comparable' to that of the domain.
\par Yet, the heart and the BZ reaction retain a partly Euclidean local
structure. At distances `far' from the inhomogeneities, is their effect truly
felt? If the anisotropy ratio is such that the `preferred' direction is only
slightly so `preferred', are spiral dynamics really affected? Can the
Oregonator distinguish the boundary from infinity if it is `very distant' from
the spiral core? \par The Euclidean model alone cannot explain these events:
clearly, any model hoping to do so should incorporate \textit{forced Euclidean
symmetry-breaking} (FESB) in order to maintain the `partly Euclidean structure'
described above. The combination of Barkley's approach with FESB predicts,
amongst other, the following (observed) spiral behaviours:\begin{description}
\item[Spiral anchoring]
appears when local inhomogeneities are present: spirals are attracted or
repelled by a RW which rotates around the site of the inhomogeneity. This has
been observed in cardiac tissue \cite{Detal} and in numerical simulations of a
modified Oregonator \cite{MPMPV}.
\item[Epicyclic drifting]
can be observed when the sizes of the physical domain and of the spiral core
are comparable: the latter is then attracted to the boundary of the domain and
rotate around it in a meandering fashion. This has been observed in experiments
and numerical simulations in a light-sensitive BZ reaction
\cite{YP,ZM}.\footnote{The term
\textit{boundary drifting} is also used.}
\item[Quasi-periodic anchoring] is
witnessed in periodically-forced RDS. The results are
similar to spiral anchoring, with the attracting/repelling structure consisting
of either two- or three-frequency quasi-periodic motion.\footnote{They are then
called \textit{entrainment} and
\textit{resonance attractors}, respectively.} These have been observed in a
light-sensitive BZ reaction which is periodically hit by light pulses and the
corresponding modified Oregonator model \cite{GZM,BSE}.
\item[Discrete RW and MRW] can be seen in systems that incorporate the notion of
anisotropy (\textit{i.e.} the system has a `preferred
direction').\footnote{Generically, these waves cannot occur in
$\GSE(2)-$equivariant reaction-diffusion systems \cite{BKT}.} In general, the
tip path has discrete two-fold symmetry. These have been observed in numerical
experiments on the bidomain model of cardiac electrophysiology~\cite{R1}.
\item[Phase-locking]
can also be seen in systems that incorporate anisotropy: the rotation and
meander frequencies can lock and this motion can be superimposed with a slow
drift. This has been observed in the bidomain model as well \cite{R1}.
\end{description}
\begin{figure}[t]
\begin{center}\includegraphics[height=210pt]{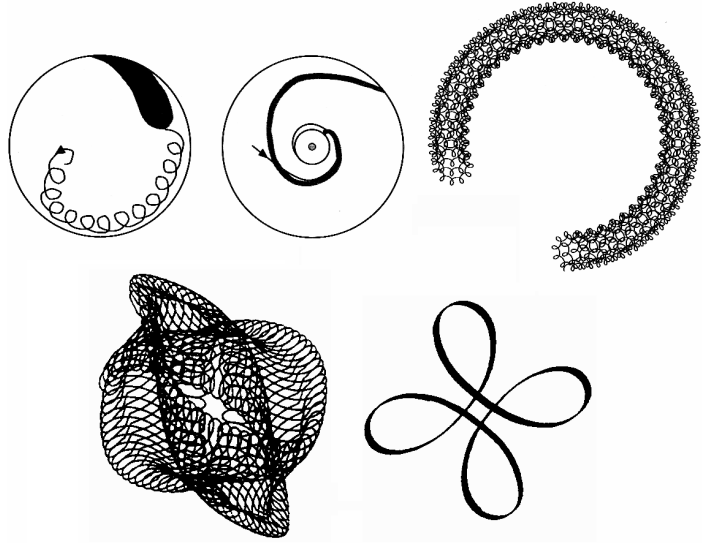}
\caption[Non-standard motions of the spiral wave tip.]{\small Non-standard
motions of the spiral wave tip.~Top row (from left to right): boundary drifting, spiral anchoring, and a three-frequency attractor \cite{LW}.\label{LW222}
  Bottom row (from left to right): a discrete MRW with $\GZ_2-$symmetry and phase-locking with small drift \cite{R1}.\label{thats}}\end{center}\hrule
\end{figure} These are illustrated in figure~\thefigure. In what follows, we
present the results obtained through systematic breaking of the Euclidean
symmetry.
\subsection{A Single TSB Term} Using the center manifold reduction theorem of
 \cite{SSW1,SSW2}, LeBlanc and Wulff
\cite{LW} showed that translational symmetry-breaking (TSB) from Euclidean
symmetry generically leads to anchoring or quasi-periodic anchoring and that
TW, MTW, boundary drifting and quasi-periodic attractors could also occur. They
did so by studying a general perturbed ODE system on the center bundle
$\GSE(2)\times \GC$, taking the form
\begin{align}\label{thefifth}
\begin{split}
\dot{p}&=e^{i\varphi}\left[F^p(q,\overline{q})+\varepsilon G^{p}(pe^{-i\varphi},\overline{p}e^{i\varphi},q,\overline{q},\varepsilon)\right]\\
\dot{\varphi}&=F^{\varphi}(q,\overline{q})+\varepsilon G^{\varphi}(pe^{-i\varphi},\overline{p}e^{i\varphi},q,\overline{q},\varepsilon) \\
\dot{q}&=F^q(q,\overline{q})+\varepsilon
G^{q}(pe^{-i\varphi},\overline{p}e^{i\varphi},q,\overline{q},\varepsilon)
\end{split}
\end{align} where $\varepsilon\in \GR$ is small and the $G-$perturbations are bounded and
uniformly continuous in $p$ and $q$. When $\varepsilon=0$, (\ref{thefifth}) is
$\GSE(2)-$equivariant under the action of (\ref{thefourth}), but for
$\varepsilon\neq 0$, the system is generally only $\GSO(2)-$equivariant: the
translational symmetry of the model has been broken. The following results are
proved (directly or in equivalent forms) in \cite{LW}.
\subsubsection{Relative Equilibria} In the case of normally hyperbolic (rotating) relative
equilibria,\footnote{That is, $q=0$ in (\ref{thefifth}) and the RW $u_*$ is not
at the transition to a MRW.} we may assume without loss of generality, and
after an appropriate time-rescaling of the $\varphi$ variable, that the center
bundle equations (\ref{thefifth}) take the form
\begin{align}\label{basiceqs4}
 \dot{p}&={\displaystyle e^{it}\left[v+ \varepsilon
H(pe^{-it},\overline{p}e^{it}, \varepsilon)\right]},
\end{align}
where $v\in \GC^{\times}$ and $\varepsilon\in\mathbb{R}$ is small.  Set
 $\widetilde{H}(w,\overline{w},\varepsilon)=
H(w-iv,\overline{w}+i\overline{v},\varepsilon)$.
\begin{thm}\label{genericthm}
Let $a=\real[D_1\widetilde{H}(0,0,0)]$.  If $a\neq 0$, then for all
$\varepsilon\neq 0$ small enough, the center bundle equation
$(\ref{basiceqs4})$ has a unique smooth branch of periodic solutions
\begin{equation}
p_{\varepsilon}(t)=\left(-iv+O(\varepsilon)\right)e^{it},\quad\varphi(t)=t,
\label{anchored_sol}
\end{equation} whose stability is exactly determined by the sign of $a\varepsilon$.
\end{thm}
These periodic solutions are centered around the origin in the $p$-plane and
are observable as anchored RW in the physical space. Note that the hypotheses
of theorem~\ref{genericthm} are generic.
\begin{thm} \label{nongenericthm} Let
$$I(\rho)=\real\left[\int_0^{2\pi}\!\!e^{-it}\widetilde{H}\left(\rho
    e^{-it},\rho e^{it},0\right)
  dt\right].$$ If $\rho_0>0$ is a hyperbolic solution of $I(\rho)=0$, then for all $\varepsilon\neq 0$ small enough,
  the center bundle equation $(\ref{basiceqs4})$ has an epicyclic solution around the origin, whose stability is
  exactly determined by the sign of $\varepsilon I'(\rho_0)$.
\end{thm}
These solutions represent quasi-periodic motion around the origin in the
$p-$plane and are observable as epicycle-like motion along a circular boundary
in the physical space, with angular frequency given by $1+O(\varepsilon)$. Note
that the hypotheses of theorem~\ref{nongenericthm} are not generic. \newl Since
TW can also be seen as $\infty-$centered orbits with no rotational component,
an appropriate change of variables considered with the limiting case $\omr\to
0$ takes the center bundle equations (\ref{thefifth}) to the equivalent system
\begin{align}\label{LWreleq3}
\dot{z}&=-vz^2+i\varepsilon z C^{\varphi}(z,\ov{z},\varepsilon)-\varepsilon z^2
C^{p}(z,\ov{z},\varepsilon),
\end{align} where $C^{\varphi}(z,\ov{z},\varepsilon)=G^{\varphi}(z^{-1},\ov{z}^{-1},\varepsilon)$ and
$C^{p}(z,\ov{z},\varepsilon)=G^{p}(z^{-1},\ov{z}^{-1},\varepsilon)$.
\begin{thm} If $C^{\varphi}$ and $C^p$ are sufficiently smooth near $z=0$, the center bundle equation $(\ref{LWreleq3})$
undergoes a transcritical bifurcation of equilibria at $\varepsilon=0$.
\end{thm}
The trivial equilibria $z=0$ represent traveling waves in the $p-$plane and are
observable as linear drifts in the physical space.
\subsubsection{Relative Periodic Solutions} In the case of normally hyperbolic (rotating) relative
periodic solutions,\footnote{That is, the $q$ equation in (\ref{thefifth}) has
a $2\pi-$periodic solution and the MRW $u^*$ is not at the transition to a RW.}
we may assume without loss of generality that, after an appropriate
time-rescaling and a change of variables, the center bundle equations
(\ref{thefifth}) take the form
\begin{align}\label{LWrelps2}
\begin{split}
\dot{w}&=-i\omr w+\varepsilon H^{w}(w,\ov{w},t,\varepsilon)\\
\dot{\varphi}&=\omr+\varepsilon H^{\varphi}(w,\ov{w},t,\varepsilon),
\end{split}
\end{align}
where $H^w$ and $H^{\varphi}$ are $2\pi-$periodic in $t$, $\omr\not\in \GZ$ and
$w=pe^{-i\varphi}$.
\begin{thm}\label{notsure}
Let $h_1^w(t)=D_wH^w(0,0,t,0)$ and set
$$\beta=\real\left[\int_{0}^{2\pi}\!\!h_1^w(t)dt\right].$$ If $\beta\neq 0$, then for all $\varepsilon\neq 0$ small
enough, the time$-2\pi$ map of the center bundle equations $(\ref{LWrelps2})$
has a unique smooth branch of hyperbolic fixed points $w_{\varepsilon}$ whose
stability is exactly determined by the sign of $\varepsilon \beta$.\end{thm}
\newpage\noindent These fixed points represent periodic solutions centered around the origin
in the $w$-plane and are observable as anchored MRW in the physical space. Note
that the hypotheses of theorem~\ref{notsure} are generic. \par There are
certain similarities between theorems~\ref{genericthm} and \ref{notsure}; in
the same vein, the following two results are related to theorem
\ref{nongenericthm}.
\begin{thm}\label{notsure2} If $\omr\not\in\GQ$, let
$$J(\rho)=\real\left[\lim_{T\to\infty}\frac{1}{T}\int_{0}^{T}\!\!e^{i\omr t}H^w\left(\rho e^{-i\omr t},\rho e^{i\omr t},t,0\right)\, dt\right].
$$ If $\rho_0>0$ is a hyperbolic solution of $J(\rho)=0,$ then for all $\varepsilon\neq 0$ small enough, the
$w-$equation in $(\ref{LWrelps2})$ has a unique smooth branch of hyperbolic
invariant two-torii $T_{\varepsilon}$, whose stability is exactly determined by
the sign of $\varepsilon J'(\rho_0)$.
\end{thm} Such an invariant two-torus represents an $O(\varepsilon)$ drift of a MRW
(centered at a point different from the origin) around its $\GSO(2)-$orbit
about $0$ in (\ref{LWrelps2}) and is observable as a three-frequency motion in
the physical space. Note that the hypotheses of theorem~\ref{notsure2} are not
generic.
\begin{thm}\label{notsure3} If $\omr \in \GQ$, with $\omr=\frac{q}{{\jmath^*}}$,
$\gcd(q,{\jmath^*})=1$ and ${\jmath^*}>1$, let
$$h_0(\xi,\ov{\xi})=\lim_{T\to\infty}\frac{1}{T}\int_{0}^{T}\!\!e^{i\omr t}H^w\left(\xi e^{-i\omr t},\ov{\xi} e^{i\omr t},t,0\right)\, dt.$$
If $\tilde{\xi}(t)$ is a hyperbolic periodic solution $[$resp. an equilibrium
point $]$ of the $\GZ_{{\jmath^*}}-$equi\-va\-riant ODE
$\dot{\xi}=h_0(\xi,\ov{\xi})$,
 then for all $\varepsilon\neq 0$ small enough, the $w-$equation in $(\ref{LWrelps2})$ has a unique smooth
branch of hyperbolic invariant two-torii $\tilde{T}(\varepsilon)$ $[$resp. of
hyperbolic $\frac{2\pi}{{\jmath^*}}-$periodic solutions
$\tilde{w}_{\varepsilon}(t)$$]$, whose stability as an invariant set is exactly
determined by the product of the sign for the stability of $\tilde{\xi}(t)$
with the sign of $\varepsilon$.\end{thm} The interpretation of these invariant
two-torii and/or periodic solutions are exactly as in the remarks following
 theorem~\ref{notsure} and \ref{notsure2}.
\newl MTW are approached much the same way as TW were tackled in the preceding section. After
an appropriate change of variables and a time-rescaling, and considering the
limiting case $\omr\to 0$, the relevant equations are
\begin{align}\label{LWmodTW3}
\begin{split}
\dot{z}&=-vz^2+\varepsilon E^z(z,\ov{z},t,\varepsilon)\\
\dot{\varphi}&=\varepsilon E^{\varphi}(z,\ov{z},t,\varepsilon),
\end{split}
\end{align} where $E^z(z,\ov{z},t,\varepsilon)=-z^2H^w(z^{-1},\ov{z}^{-1},t,\varepsilon)$,
$E^{\varphi}(z,\ov{z},t,\varepsilon)=H^{\varphi}(z^{-1},\ov{z}^{-1},t,\varepsilon)$
are $2\pi-$periodic in $t$ and $v\in \GC^{\times}$.
\begin{thm} If $E^z$ and $E^{\varphi}$ are sufficiently smooth near $z=0$, the time$-2\pi$ map of the $z-$equation
in the center bundle equations $(\ref{LWmodTW3})$ undergoes a transcritical
bifurcation of fixed points at $\varepsilon=0$. \end{thm} The trivial fixed
points $z=0$ represent MTW in the $p-$plane and are observable as linear
meandering in the physical space.

\subsection{A Single RSB Term} LeBlanc \cite{LeB} then showed that rotational symmetry-breaking (RSB) from
Euclidean symmetry could provide an explanation for the appearance of discrete
RW, discrete MRW and phase-locking in excitable media with anisotropy. He did
so by studying a general perturbed ODE system on the center bundle
$\GSE(2)\times \GC$, taking the form
\begin{align}\label{thesixth}
\begin{split}
\dot{p}&=e^{i\varphi}\left[F^p(q,\overline{q})+\varepsilon G^{p}(\varphi,q,\overline{q},\varepsilon)\right]\\
\dot{\varphi}&=F^{\varphi}(q,\overline{q})+\varepsilon G^{\varphi}(\varphi,q,\overline{q},\varepsilon) \\
\dot{q}&=F^q(q,\overline{q})+\varepsilon
G^{q}(\varphi,q,\overline{q},\varepsilon)
\end{split}
\end{align} where $\varepsilon\in \GR$ is small, ${\jmath^*}\in \GN$, and the $G-$perturbations are bounded
and uniformly continuous in $p$ and $q$, as well as
$\frac{2\pi}{{\jmath^*}}-$periodic in $\varphi$. When $\varepsilon=0$,
(\ref{thesixth}) is $\GSE(2)-$equivariant under the action of
(\ref{thefourth}), but for $\varepsilon\neq 0$, the system is generally only
$\GC\sdp \GZ_{{\jmath^*}}-$e\-qui\-va\-riant: the rotational symmetry of the
model has been broken. The following results are proved (directly or in
equivalent forms) in \cite{LeB}.
\subsubsection{Relative Equilibria} Let ${\jmath^*}\geq 1$ be a fixed integer. In the case of normally hyperbolic (rotating) relative
equilibria, we may assume without loss of generality that, after an appropriate
time-rescaling of the $\varphi$ variable, the center bundle equations
(\ref{thesixth}) take the form
\begin{align}\label{Lreleq}
\begin{split}
\dot{p}&=e^{i\varphi}[v+\varepsilon G^{p}(\varphi,\varepsilon)]\\
\dot{\varphi}&=\omr+\varepsilon G^{\varphi}(\varphi,\varepsilon),
\end{split}
\end{align} where $v\in \GC^{\times}$, $\omr\in \GR$ and $G^p,G^{\varphi}$ are $\frac{2\pi}{{\jmath^*}}-$periodic
 in $\varphi$. \begin{thm}
 Assume $\omr\neq 0$. For all $\varepsilon\neq 0$ sufficiently small, the solutions of $(\ref{Lreleq})$ are
 $\frac{2\pi}{{\jmath^*}}-$periodic in time with discrete $\GZ_{{\jmath^*}}-$symmetry.
 \end{thm} These periodic solutions represent discrete RW in the physical
space.
\begin{thm} Assume $\omr=0$. For all $\varepsilon\neq 0$ sufficiently
small, if $G^{\varphi}(\varphi,0)\neq 0$ for all $\varphi\in \GS^1$, the
solutions of $(\ref{Lreleq})$ are discrete RW, with (large) radii of the order
of $\frac{1}{\varepsilon}$. On the other hand, if there exists $\varphi^*\in
[0,2\pi)$ such that $G^{\varphi}(\varphi^*,0)=0$ and
$D_{\varphi}G^{\varphi}(\varphi^*,0)\neq 0$, then $(\ref{Lreleq})$ has at least
${\jmath^*}$ stable (attracting) TW solutions and an equal number of unstable
(repelling) TW solutions.
\end{thm}
In the latter case, all solutions of (\ref{Lreleq}) end up drifting linearly,
after an initial transient period.
\subsubsection{Relative Periodic Solutions}
In the case of normally hyperbolic (rotating) relative periodic
solutions,footnote{Where the corresponding $2\pi-$periodic solution $q^*(t)$ to
the $q$ equation in (\ref{thesixth}) is such that
$F^q(q^*(t),\overline{q^*}(t))\neq 0$ for all $t$.} we may assume without loss
of generality that, after an appropriate time-rescaling and a change of
variables, the center bundle equations (\ref{thesixth}) take the form
\begin{align}\label{Lrelps2}
\begin{split}
\dot{p}&=e^{i\varphi}[\tilde{F}^p(\theta)+\varepsilon \tilde{G}^{p}(\varphi,t,\varepsilon)]\\
\dot{\varphi}&=\omr+\tilde{F}^{\varphi}(\theta)+\varepsilon
\tilde{G}^{\varphi}(\varphi,t,\varepsilon),
\end{split}
\end{align}
where all functions are $2\pi-$periodic in $\theta$ and
$\frac{2\pi}{{\jmath^*}}-$periodic in $\varphi$, and where the average value of
$\tilde{F}^{\varphi}$ is $0$. In the unperturbed case, all solutions of
(\ref{Lrelps2}) are MTW (if $\omr\in \GZ$) or discrete MRW (if $\omr\not\in
\GZ$). \newpage\noindent If $\varepsilon\neq 0$ is sufficiently small, the
$\varphi-$equation in (\ref{Lrelps2}) defines a $\GZ_{{\jmath^*}}-$equivariant
flow on a two-torus, with associated Poincar\'e map
\begin{equation}\label{Pmap2T} P(\varphi;\omr,\varepsilon)=\varphi+2\pi\omr+\varepsilon H(\varphi,\omr,\varepsilon),\end{equation}
where $H$ is $\frac{2\pi}{{\jmath^*}}-$periodic in $\varphi$ and $C^r$ for some
$r\geq 3$. The map $P$ is thus a circle map; denote its rotation number by
$\rho_{\omr,\varepsilon}$.
\begin{thm} If
$\rho_{\omr,\varepsilon}=\frac{m}{\gamma}$, where $\gamma>0,m\in \GZ$ are
coprime, set $k=\gcd(\gamma,{\jmath^*})$. Then, phase-locking occurs in
$(\ref{Lrelps2})$ when $(\omr,\varepsilon)$ lies in the $m:\gamma$ Arnol'd
tongue of $P$.\end{thm} If $k\neq 1$, the $p-$component of solutions of
$(\ref{Lrelps2})$ is a $\gamma-$petaled $2\pi \gamma-$periodic curve with
$\GZ_k-$spatial symmetry. Otherwise, the $p-$component of solutions of
$(\ref{Lrelps2})$ is a superposition of a motion akin to the one in the case
$k\neq 1$ together with a `slow' linear drift.
\begin{thm} If $\rho_{\omr,\varepsilon}\not\in \GQ$, the flow on the
two-torus described above is ergodic.\end{thm}
If there exist $\sigma\in (0,1)$
and $K>0$ such that
$$\left|\rho_{\omr,\varepsilon}+\frac{k}{j}\right|\geq
K|j|^{-(2+\sigma)}$$ for all non-trivial integer pairs $(k,j)$ (which is almost
always the case), then the $p-$component of solutions of (\ref{Lrelps2}) is
quasi-periodic and the closure of its positive image has
$\GZ_{{\jmath^*}}-$rotational symmetry.
\subsection{Simultaneous TSB Terms} The next logical step lies in studying the effects of $n$ simultaneous TSB
perturbations, for $n>1$, which is done in \cite{BLM,Bo2,Byeah}. Any excitable
media which is littered with inhomogeneities, such as the human heart, could
then in theory be modeled by a general system of ODE on the center bundle
$\GSE(2)\times \GC$, taking the form
\begin{align}\label{theseventh}
\begin{split}
\dot{p}&=e^{i\varphi}\left[F^p(q,\overline{q})+\sum_{i=1}^n\lambda_i G^{p}_i((p-\xi_i)e^{-i\varphi},\overline{(p-\xi_i)}e^{i\varphi},q,\overline{q},\lambda)\right]\\
\dot{\varphi}&=F^{\varphi}(q,\overline{q})+\sum_{i=1}^n\lambda_i G^{\varphi}_i((p-\xi_i)e^{-i\varphi},\overline{(p-\xi_i)}e^{i\varphi},q,\overline{q},\lambda) \\
\dot{q}&=F^q(q,\overline{q})+\sum_{i=1}^n\lambda_i
G^{q}_i((p-\xi_i)e^{-i\varphi},\overline{(p-\xi_i)}e^{i\varphi},q,\overline{q},\lambda)
\end{split}
\end{align} where $\lambda=(\lambda_1,\ldots,\lambda_n)\in \GR^n$ is small, $\xi_1,\ldots,\xi_n\in \GC$ are all distinct and the $G-$perturbations are bounded and
uniformly continuous in $p$ and $q$. \par When $\lambda=0$, (\ref{theseventh})
is $\mathbb{S}\mathbb{E}(2)$-equivariant under the action of (\ref{thefourth});
when $\lambda\neq 0$ is near the origin and along the $j^{\mbox{\footnotesize
th}}$ coordinate axis of $\mathbb{R}^n$,  (\ref{theseventh}) is generally only
$\mathbb{S}\mathbb{O}(2)_{\xi_j}$-equivariant (i.e. it only commutes with
rotations about the point $\xi_j$), and when two or more of the $\lambda_i$ are
not zero, the system is generally only trivially equivariant: the translational
symmetry of the model has been broken.
\subsubsection{Relative Equilibria}\label{re}
In the case of normally hyperbolic (rotating) relative equilibria, we may
assume without loss of generality that, after an appropriate time-rescaling of
the $\varphi$ variable, the center bundle equations (\ref{theseventh}) take the
form
\begin{align} \label{system1}
\dot{p}&=e^{it}\left[v+
\sum_{j=1}^n\,\lambda_jH_j((p-\xi_j)e^{-it},(\overline{p-\xi_j})e^{it},
\lambda)\right]
\end{align}
where, $v\in \GC^{\times}$ and the functions $H_j$ are smooth and uniformly
bounded in $p$. Boily, LeBlanc and Matsui study spiral anchoring in this
particular setting \cite{BLM,Byeah}.
\par
A $2\pi-$periodic solution $p_{\lambda}$ of $(\ref{system1})$ is called a
\textit{perturbed rotating wave} of $(\ref{system1})$. Define the average value
\begin{align*}[p_{\lambda}]_{\A}&=\frac{1}{2\pi}\int_{0}^{2\pi}\!\!\!\!p_{\lambda}(t)\, dt.\end{align*}
If the Floquet multipliers of $p_{\lambda}$ all lie within (resp. outside) the
unit circle, we shall say that $[p_{\lambda}]_{\A}$ is the \textit{anchoring}
(resp. \textit{repelling}, or \textit{unstable anchoring}) \textit{center}
of~$p_{\lambda}$.
\begin{thm} \label{genericthm2} Let $k\in \{1,\ldots,n\}$ and define
$\alpha_k=\real\left[D_1H_k(iv,-i\ov{v},0)\right]$. If $\alpha_k\neq 0$, there
exists a wedge-shaped region of the form
\[
{\mathcal W}_{k}=\{(\lambda_1,\ldots,\lambda_n)\in
\GR^n\,:\,|\lambda_j|<W_{k,j}|\lambda_k|,\,\,\,W_{k,j}>0,\,\,\mbox{\rm for
$j\neq k$ and $\lambda_k$ near}\,\,0\,\}
\]
such that for all $0\neq \lambda\in {\mathcal W}_{k}$, the center bundle
equation $(\ref{system1})$ has a unique perturbed rotating wave
$\mathcal{S}^k_{\lambda}$, with center $[\mathcal{S}^k_{\lambda}]_{\A}$
generically away from $\xi_k$, whose stability is uniquely determined by the
sign of $\lambda_k\alpha_k$.
\end{thm}
In contrast to theorem~\ref{genericthm}, these periodic solutions are not
necessarily centered around an inhomogeneity $\xi_k$ in the $p-$plane, but they
 are still observable as anchored RW in the physical space. Note that the hypotheses
of theorem~\ref{genericthm2} are generic. Furthermore,
$[\mathcal{S}^k_{\lambda}]_{\A}$ is a center of anchoring when
$\lambda_k\alpha_k< 0$ and a center of repelling when $\lambda_k\alpha_k> 0$.
\par Boily then showed that theorem~\ref{nongenericthm} has a similar
generalization \cite{Bo2,Byeah}. An \textit{epicycle manifold} of
(\ref{system1}) is an invariant set~$\hat{\mathcal{E}}_{\lambda}$ for
(\ref{system1}) in which all solutions are epicycles when projected upon the
$p-$plane.
\begin{thm} \label{nongenericthm2} Let $k\in \{1,\ldots, n\}$ and $$I_k(\rho)=\real\left[\int_0^{2\pi}\!\!e^{-it}H_k\left(\rho
    e^{-it}\!\!-iv,\rho e^{it}\!\!+i\ov{v},0\right)
  dt\right].$$ If $\rho_k>0$ is a hyperbolic solution of $I_k(\rho)=0$, then there exists a wedge-shaped region of the form
\[
{\mathcal{V}}_{k}=\{(\lambda_1,\ldots,\lambda_n)\in
\GR^n\,:\,|\lambda_j|<V_{k,j}|\lambda_k|,\,\,\,V_{k,j}>0,\,\,\mbox{\rm for
$j\neq k$ and $\lambda_k$ near}\,\,0\,\}
\]
such that for all $0\neq \lambda\in {\mathcal{V}}_{k}$, $(\ref{system1})$ has
an epicycle manifold ${\mathcal{E}}_{\lambda}^k$ whose stability is exactly
determined by the sign of $\lambda_kI'_k(\rho_k)$.
\end{thm}
These solutions represent quasi-periodic motion centered near (but not
generically at) the inhomogeneity $\xi_k$ in the $p-$plane and are observable
as epicycle-like motion along a circular boundary in the physical space. Note
that as was previously the case, the hypotheses of theorem~\ref{nongenericthm2}
are not generic.
\subsection{Combined TSB and RSB Terms} Yet another way in which the Euclidean
symmetry can be broken lies in the combination of TSB and RSB terms; such a
situation is analyzed in \cite{Bo1,Byeah}. Anisotropic media near an
inhomogeneity, such as cardiac tissue in the neighbourhood of site of higher
density at the origin, could then in theory be modeled by a general system of
ODE on the center bundle $\GSE(2)\times \GC$, taking the form
\begin{align}\label{theeight}
\begin{split}
\dot{p}&=e^{i\varphi}\left[F^p(q,\overline{q})+\varepsilon G^p(\varphi,q,\ov{q},\varepsilon,\mu)+\mu H^p(pe^{-i\varphi},\ov{p}e^{i\varphi},q,\ov{q},\varepsilon,\mu)\right] \\
\dot{\varphi}&=F^{\varphi}(q,\overline{q})+\varepsilon
G^{\varphi}(\varphi,q,\ov{q},\varepsilon,\mu)+\mu
H^{\varphi}(pe^{-i\varphi},\ov{p}e^{i\varphi},q,\ov{q},\varepsilon,\mu)\\
\dot{q}&=F^{q}(q,\overline{q})+\varepsilon
G^{q}(\varphi,q,\ov{q},\varepsilon,\mu)+\mu
H^{q}(pe^{-i\varphi},\ov{p}e^{i\varphi},q,\ov{q},\varepsilon,\mu),
\end{split}
\end{align}
where $(\varepsilon,\mu)\in \GR^2$ is small, ${\jmath^*}\in \GN$, the
$G,H-$perturbations are bounded and uniformly continuous in $p$ and $q$, and
the $G-$perturbations are $\frac{2\pi}{{\jmath^*}}-$periodic in $\varphi$.
Throughout this section, we fix ${\jmath^*}\in \GN$. \par Let $\GC\sdp
\GZ_{{\jmath^*}}$ be the subgroup of $\GSE(2)$ containing all translations and
rotations by angle $\frac{2\pi k}{{\jmath^*}}$, $k\in \GZ$. When
$(\varepsilon,\mu)=0$, (\ref{theeight}) is
$\mathbb{S}\mathbb{E}(2)$-equivariant under the action of (\ref{thefourth});
when $\varepsilon\neq 0$ is small and $\mu=0$, it is $\GC\sdp
\GZ_{{\jmath^*}}-$equivariant; when $\mu\neq 0$ is small and $\varepsilon=0$,
it is $\GSO(2)-$equivariant, and it is generally only trivially equivariant
otherwise: both the translational symmetry and the rotational symmetry of the
model has been broken. \par In the case of normally hyperbolic (rotating)
relative equilibria, we may assume without loss of generality that, after an
appropriate time-rescaling of the $\varphi$ variable, the center bundle
equations (\ref{theseventh}) take the form
\begin{equation}\label{H*R.com}\dot{p}=e^{it}\big[v+\varepsilon G(t,\varepsilon,\mu)+\mu H(pe^{-it},\ov{p}e^{it},\varepsilon,\mu)\big],\end{equation}
 where $(\varepsilon,\mu)\in\GR^2$, $v\in \GC$, $G$ and $H$ are smooth and uniformly bounded in $p$, and $G$ is $\frac{2\pi}{\jmath^*}-$periodic in $t$.
The following results are proved (directly or in equivalent forms) in
\cite{Bo1,Byeah}; they depend greatly on the nature of $\jmath^*$.
\subsubsection{The Case ${\jmath^*}= 1$} As $G$ is then $2\pi-$periodic in $t$, it can be written as the uniformly
convergent Fourier series
\begin{align}\label{FourierG1}
G(t,\varepsilon,\mu)=\sum_{n\in {\mathbb Z}}\, g_n(\varepsilon,\mu)e^{int}.
\end{align}
 \begin{thm} \label{noclueagain} If $g_{-1}(0,0)\neq 0$ and if $c_1=\real[D_1H(-iv,i\ov{v},0,0)]\neq
 0$,
 there exists a wedge-shaped region of the form
\[
{\mathcal W}=\{(\varepsilon,\mu)\in
\GR^2\,:\,|\varepsilon|<K|\mu|,\,\,\,K>0,\,\,\mbox{$\mu$ near}\,\,0\,\}
\]
such that for all $(\varepsilon,\mu)\in {\mathcal W}$ with $\varepsilon\neq 0$,
the center bundle equation $(\ref{H*R.com})$ has a unique hyperbolic discrete
rotating wave ${\mathcal D}^{1}_{\varepsilon,\mu}$ with trivial spatio-temporal
symmetry, centered near but generically not at the origin, whose stability is
exactly determined by the sign of  $\mu c_1$.
\end{thm}
\begin{thm} \label{whathte} Let
$$R(\rho)=\real\left[\int_0^{2\pi}\!\!e^{-it}H\left(\rho
    e^{-it}\!\!-iv,\rho e^{it}\!\!+iv,0,0\right)
  dt\right].$$ If $\rho_0>0$ is a hyperbolic solution of $R(\rho)=0$, then there exists a wedge-shaped region  of the form
$${\mathcal V}=\{(\varepsilon,\mu)\in \GR^2\,:\,|\varepsilon|<K|\mu|,\,\,\,K>0,\,\,\mbox{$\mu$ near}\,\,0\,\}
$$
such that for all $(\varepsilon,\mu)\in {\mathcal V}$ with $\varepsilon\neq 0$,
the center bundle equation $(\ref{H*R.com})$ has an epicycle manifold
${\mathcal G}_{\varepsilon,\mu}^1$ centered near but generically not at the
origin, whose stability is exactly determined by the sign of $\mu R'(\rho_0)$.
\end{thm}
To get the complete picture (in both of these theorems), the situation would
also need to be analyzed near the $\varepsilon-$axis. This implies dealing with
fixed points of maps at $\infty$; as of now, it has been relegated to a future
 investigation. Note that the hypotheses of theorem~\ref{noclueagain} are
generic, while those of~\ref{whathte} are not.

\subsubsection{The Case
${\jmath^*} >1$} In fully anisotropic media, one would have $\jmath^*=2$.
 \begin{thm} \label{noclue} Let $c_1$ be as in theorem~$\ref{noclueagain}$. If $c_1\neq 0$, then there
 exists a small deleted neighbourhood
   $\mathcal{W}^{{\jmath^*}}$ of the origin, such that for all $0\neq (\varepsilon,\mu)\in \mathcal{W}^{{\jmath^*}}$, the center bundle equation
   $(\ref{H*R.com})$ has a unique hyperbolic discrete rotating wave ${\mathcal D}^{{\jmath^*}}_{\varepsilon,\mu}$
   with $\GZ_{{\jmath^*}}-$spatio-temporal symmetry, centered at the origin in the $p-$plane, whose stability
   is exactly determined by the sign of $\mu c_1$.
\end{thm}
\begin{thm}\label{ahahahthmtorus} Let $R(\rho)$ be as in
theorem~$\ref{whathte}$. If $\rho_0>0$ is a hyperbolic solution of $R(\rho)=0$,
then there exists a small deleted neighbourhood $\mathcal{V}^{{\jmath^*}}$ of
the origin, such that for all $(\varepsilon,\mu)\in \mathcal{V}^{{\jmath^*}}$,
the center bundle equation $(\ref{H*R.com})$ has an epicycle manifold
${\mathcal G}^{{\jmath^*}}_{\varepsilon,\mu}$ centered at the origin, whose
stability is exactly determined by the sign of $\mu c_1$.
\end{thm} Contrary to what might be thought at first, the epicycle manifolds of theorem \ref{ahahahthmtorus} do \textsl{not} have
$\GZ_{{\jmath^*}}-$spatio-temporal symmetry; however, the epicycles themselves
possess this symmetry in an appropriate frame of reference. As has been the
case throughout this article, the hypotheses of theorem~\ref{noclue} are
generic, while those of theorem~\ref{ahahahthmtorus} are not.

%%% ----------------------------------------------------------------------

\section{Conjectures, Related Results and Future Work}
The application of FESB to the study of rotating waves has yielded a number of
interesting verifiable results: the predictive power of the dynamical approach
cannot be denied. Yet the full picture of spiral wave dynamics is far from
complete
\subsection{Conjectures}\label{C&R} In this section, we present some conjectures concerning modulated rotating
waves.
\subsubsection{Simultaneous TSB Terms} The equations describing the essential dynamics of a normally
hyperbolic modulated rotating wave are similar to those of rotating waves. Near
such a MRW, the center bundle equations (\ref{thesixth}) are equivalent to
\begin{align} \label{system12conj}
\begin{split}
\dot{p}&=e^{i\varphi}\left[v+\sum_{j=1}^n\lambda_jH_j^p((p-\xi_j )e^{-i\varphi},\overline{(p-\xi_j)}e^{i\varphi},t,\lambda)\right] \\
\dot{\varphi}&=\omr+\sum_{j=1}^n\lambda_jH_j^{\varphi}((p-\xi_j
)e^{-i\varphi},\overline{(p-\xi_j)}e^{i\varphi},t,\lambda),
\end{split}
\end{align} where $\lambda=(\lambda_1,\ldots,\lambda_n)\in \GR^n$, $v\in \GC$, $\omr\neq 0$ and the functions
$H_j^{p,\varphi}$  are smooth, uniformly bounded in $p$ for all $j=1,\ldots, n$
and $2\pi-$periodic in $t$.\par System (\ref{system12conj}) cannot be analyzed
as easily as (\ref{system1}), but it seems nonetheless likely that the
following hold.
\begin{conj} \label{thm43conj} Let $k\in \{1,\ldots,n\}$. Generically,
 there is a wedge-shaped region of the form
\[
{\mathcal W}_{k}=\{(\lambda_1,\ldots,\lambda_n)\in
\GR^n\,:\,|\lambda_j|<W_{k,j}|\lambda_k|,\,\,\,W_{k,j}>0,\,\,\mbox{\rm for
$j\neq k$ and $\lambda_k$ near}\,\,0\,\}
\]
such that for all $0\neq \lambda\in {\mathcal W}_{k}$, the center bundle
equations $(\ref{system12conj})$ have a unique perturbed modulated rotating
wave solution $\mathfrak{S}^k_{\lambda}$, centered generically away from
$\xi_k$.
\end{conj}
\begin{conj}\label{thmnld2conj} Let $k\in \{1,\ldots, n\}$. Given a hyperbolic equilibrium of a related averaged system,
there is a wedge-shaped region  of the form
\[
{\mathcal V}_{k}=\{(\lambda_1,\ldots,\lambda_n)\in
\GR^n\,:\,|\lambda_j|<V_{k,j}|\lambda_k|,\,\,\,V_{k,j}>0,\,\,\mbox{\rm for
$j\neq k$ and $\lambda_k$ near}\,\,0\,\}
\]
such that for all $0\neq \lambda\in {\mathcal V}_{k}$, the center bundle
equations $(\ref{system12conj})$ have a hyperbolic $3-$frequency epicycle
manifold ${\mathfrak{E}}_{\lambda}^k$ centered near, but generically not at,
$\xi_k$.\end{conj}
\subsubsection{Combined RSB and TSB Terms}
Fix $\jmath^*\in \GN$. The equations describing the essential dynamics of a
normally hyperbolic modulated rotating wave are similar to those of rotating
waves. Near such a MRW, the center bundle equations (\ref{thesixth}) are
equivalent to
\begin{align} \label{system12conj1}
\begin{split}
\dot{p}&=e^{i\varphi}\left[v+\varepsilon G^p(\varphi,t,\varepsilon,\mu)+\mu H^p(pe^{-i\varphi},\ov{p}e^{i\varphi},t,\varepsilon,\mu)\right] \\
\dot{\varphi}&=\omr+\varepsilon G^{\varphi}(\varphi,t,\varepsilon,\mu)+\mu
H^{\varphi}(pe^{-i\varphi},\ov{p}e^{i\varphi},t,\varepsilon,\mu),
\end{split}
\end{align} where $(\varepsilon,\mu)\in\GR^2$, $v\in \GC$, $\omr\neq 0$, the functions $G^{p,\varphi},H^{p,\varphi}$ are smooth and
uniformly bounded in $p$ and $2\pi-$periodic in $t$, and the functions
$G^{p,\varphi}$ are additionally $\frac{2\pi}{\jmath^*}-$periodic in $\varphi$.
\par System (\ref{system12conj1}) cannot be analyzed as easily as (\ref{H*R.com}),
it seems nonetheless likely that the following  hold.
 \begin{conj} \label{noclueagainconj}Let ${\jmath^*}=1$. Generically, there is a wedge-shaped region of the form
\[
{\mathcal W}=\{(\varepsilon,\mu)\in
\GR^2\,:\,|\varepsilon|<K|\mu|,\,\,\,K>0,\,\,\mbox{$\mu$ near}\,\,0\,\}
\]
such that for all $(\varepsilon,\mu)\in {\mathcal W}$ with $\varepsilon\neq 0$,
the center bundle equations $(\ref{system12conj1})$ have a unique hyperbolic
discrete modulated rotating wave $\mathfrak{D}^{1}_{\varepsilon,\mu}$, with
trivial spatio-temporal symmetry, centered away from the origin.
\end{conj}
\begin{conj}\label{ahahthmtorusconj} Let ${\jmath^*}=1$. Given a hyperbolic equilibrium of a related averaged system,
there is a wedge-shaped region  of the form
$${\mathcal V}=\{(\varepsilon,\mu)\in \GR^2\,:\,|\varepsilon|<K|\mu|,\,\,\,K>0,\,\,\mbox{$\mu$ near}\,\,0\,\}
$$
such that for all $(\varepsilon,\mu)\in {\mathcal V}$ with $\varepsilon\neq 0$,
the center bundle equations $(\ref{system12conj1})$ have a hyperbolic
$3-$frequency epicycle manifold $\mathfrak{G}_{\varepsilon,\mu}^1$ centered
near but generically not at the origin.
\end{conj}
The remark following theorem \ref{whathte} is likely to hold for
theorems~\ref{noclueagainconj} and \ref{ahahthmtorusconj}.
 \begin{conj} \label{noclueconj}Let ${\jmath^*}>1$. Generically, there is a deleted neighbourhood $\mathcal{W}^{{\jmath^*}}$ of the origin such that
 such that for all $(\varepsilon,\mu)\in {\mathcal W}^{{\jmath^*}}$, $(\ref{system12conj1})$ has a unique hyperbolic discrete modulated rotating wave
 $\mathfrak{D}^{{\jmath^*}}_{\varepsilon,\mu}$ with $\GZ_{{\jmath^*}}-$spatio-temporal symmetry centered at the origin.
\end{conj}
\begin{conj}\label{ahahahthmtorusconj} Let ${\jmath^*}>1$. Given a hyperbolic equilibrium of a related averaged system,
there is a deleted neighbourhood $\mathcal{V}^{{\jmath^*}}$ of the origin such
that for all $(\varepsilon,\mu)\in {\mathcal V}^{{\jmath^*}}$,
$(\ref{system12conj1})$ has a hyperbolic $3-$frequency epicycle manifold
$\mathfrak{G}^{{\jmath^*}}_{\varepsilon,\mu}$ centered at the origin.
\end{conj}
\subsection{Remarks} While numerical experiments show without the shadow of a doubt
that spirals can anchor away from a center of inhomogeneity
\cite{BLM,Bo1,Byeah}, this writer would find it very satisfying to see this
result reproduced in the laboratory. From a resolutely profane perspective, the
Belousov-Zhabotinsky reaction appears most likely to yield results, as it seems
the easiest to control.\par It should also be noted that the study of spirals
does not start and end in the plane. For instance, Comanici used the dynamical
system approach in her doctoral thesis to study spirals on spherical domains
\cite{C1}. Scroll waves, the 3-dimensional analogues of spiral waves, have also
attracted attention from the physics and cardiology communities in recent years
\cite{Scroll1,Scroll2,Scroll3,Scroll4}.
\subsection{Future and Related Work} We finish this paper with a list of problems
and open questions: their solutions would improve our knowledge and
understanding of spiral waves and excitable media.
\subsubsection{MRW, TW and MTW}
When the FESB becomes too complex, only RW and TW are easily amenable to
characterization. Even then, TW have not yet been tackled. The current approach
\cite{BLM,Bo1,Bo2} needs to be suitably modified so as to accommodate the extra
variable that appears in (\ref{system12conj}) and (\ref{system12conj1}).
\subsubsection{An Explicit Center Manifold Reduction Theorem}
While extremely powerful, the CMRT has the disadvantage of being a strict
existence theorem: it tells us that the dynamics on the center manifold are
given by an ODE system with certain symmetries, but it does not provide the
explicit relation between that system and the original semi-flow.
\par In particular, when a specific center manifold system is studied, we have
no way of knowing if it corresponds to a `viable' excitable system,
\textsl{i.e.} if it is `attainable' in any way from such an excitable system
\emph{via} the CMRT. \par Recent observations by Lajoie and LeBlanc \cite{LL}
suggest that it might be possible to efficiently relate the coefficients of a
RDS to those of the center manifold near a traveling wave. Is there a direct
and efficient way to compute the relevant coefficients of the center bundle
equations directly from the PDE, near any type of relative equilibrium or
relative periodic solution?
\subsubsection{Spiral Groupings}
So far, the model-independent approach based on the CMRT has only been used to
study \emph{isolated} spiral waves. Even though experiments by Li, Ouyang,
Petrov and Swinney \cite{LOPS} have shown that spiral waves can be isolated
with the help of a laser, they are rarely found in that state in excitable
media (see for instance the two illustrations on p.~\pageref{Ballill}). \par
Spiral groupings, where two or more spirals rotate around a common center or
one another, have much different dynamics, as can be attested by the recent
numerical simulations of Pertsov and Zariski \cite{PZ}. Some of the showcased
interactions are somewhat analogous to already-obtained results about epicycle
drifting, which begs the question: how can the current approach be altered to
apply to spiral groupings as well?
\subsubsection{Global Spiral Dynamics}
Finally, the CMRT can only be applied to local neighbourhoods of spiral wave
solutions. Yet, many spiral interactions are global in nature (see \cite{PZ}
for details). As of now, there is little machinery short of numerical
simulations to deal with spiral dynamics on a global level. How can this
situation be remedied?

\appendix

\section{Selected Bibliography of Spiral Pattern
Formation}
\begin{description}
\item[BZ Reaction and the Oregonator $-$]
\cite{LOPS,GZM,ZM,MPMPV}.
\item[Cardiac Tissue and the FHN Equations $-$]
\cite{WR,YP,Detal,R1,Wetal,J,MAK}.
\item[Global Spiral Dynamics $-$]
\cite{Hendrey,PZ}.
\item[Surveys $-$]
\cite{JW,W1}.
\item[Dynamical System Approach $-$]
\cite{BKT,B1,BK,SSW2,GLM,Wulff,SSW3,GLM2,LW,LeB,C1,BLM,Bo1,Bo2,Byeah}.
\end{description}

% ------------------------------------------------------------------------
%Included for Gather Purpose only:
\bibliographystyle{amsplain}
\bibliography{../Bibliography/spirals}
\end{document}